%% file: reduction-en.tex
\documentclass[11pt]{article}
\usepackage{amsmath,amssymb,latexsym,a4}

\reversemarginpar

\input{macros.tex}

\renewcommand{\GLP}{\mathrm{GLP}}
\renewcommand{\J}{\mathrm{J}}

\begin{document}

\title{Notes on a reduction property for GLP-algebras}

\author{L. Beklemishev\thanks{Research supported by the Russian Science Foundation under grant No.~14--50--00005.} \\[1ex]
{\small \parbox{10cm}{
\begin{center}
Steklov Mathematical Institute \\
Gubkina str. 8, 119991 Moscow, Russia \\
e-mail: \texttt{bekl@mi.ras.ru}
\end{center}}
}}

\maketitle

\begin{abstract}
We consider some natural generalizations to the class of all GLP-algebras of the so-called reduction property for reflection algebras in arithmetic. An analogue of this property is established for the free GLP-algebras and for some topological GLP-algebras (GLP-spaces).
\end{abstract}

The notion of GLP-algebra emerged in the study of modal logics describing the behaviour of reflection principles and formalized $\gw$-consistency in Peano arithmetic. From the point of view of modal logic, GLP-algebras are models of polymodal provability logic $\GLP$ introduced by G.~Japaridze \cite{Dzh86,Dzh88}. Central examples of GLP-algebras are reflection algebras associated with formal arithmetical theories $T$  (originally called \emph{graded provability algebras} in \cite{Bek04,Bek05}) described below.

\emph{$\Sigma_n$-reflection formula} for an arithmetical r.e.\ theory $T$ (containing $\EA=I\Delta_0+\exp$) is a sentence $R_n(T)$, expressing in the language of Peano arithmetic the fact that each $T$-provable sentence of arithmetical complexity $\Sigma_n$ is true. Every such formula correctly defines an operation $\la n\ra:\cL_T\to \cL_T$ on the Lindenbaum boolean algebra $\cL_T$ of $T$ which associates with the equivalence class of a sentence $\phi$ the class of $R_n(T+\phi)$:
 $$\la n\ra: [\phi] \longmapsto [R_n(T+\phi)].$$
The Lindenbaum algebra of $T$ enriched by all these operations is called the \emph{reflection algebra of $T$} and is denoted $\cM_T=(\cL_T,\{\la n\ra :n<\gw\})$.

Interesting applications of reflection algebras are connected with the proof-theoretic study of Peano arithmetic $\PA$ and, in particular, with the description of the set of theorems of $\PA$ of quantifier complexity $\Pi_n^0$ in terms of a natural system of ordinal notation for the ordinal $\ge_0=\sup\{\gw,\gw^\gw,\dots\}$.

The reflection algebra of $T$, in addition to the identities of boolean algebras, satisfies the following principles (where we denote $[n]x:= \neg \la n\ra \neg x$):
\benr
\item $\la n\ra(x\lor y)=(\la n \ra x\lor \la n \ra y)$;
\item $\la n \ra 0=0$;
\item $\la n\ra x = \la n\ra (x\land \neg \la n\ra x)$; \label{lob}
\item $\la n \ra x \leq \la m \ra x$, for $m\leq n$; \label{p1}
\item $\la m \ra x \leq [n]\la m\ra x$, for $m<n$. \label{p2}
\eenr
Notice that \emph{L\"ob's identity} \refeq{lob} formalizes a generalization of G\"odel's second incompleteness theorem for $T$.

\bd \emph{GLP-algebra} is a boolean algebra $\cM$ enriched by the operations
    $\la n \ra$, for all $n<\gw$, satisfying identities (i)--(v).
\ed
Thus, reflection algebras are examples of GLP-algebras. Other examples emerge in the study of scattered topological spaces, in particular, in the study of ordinal topologies \cite{BekGab13,BekGab14}. Free GLP-algebras present an independent interest in connection with the study of polymodal provability logic GLP.

The reduction property for reflection algebras is a key fact needed for the proof-theoretic analysis of Peano arithmetic. In \cite{Bek04} this property was formulated using the notion of $\Pi^0_n$-conservativity which pertains to formal arithmetic but has no exact correspondent on the abstract algebraic level. In this note we propose some natural analogues of the reduction  property that make sense for arbitrary GLP-algebras and establish it for some other classes of GLP-algebras, not necessarily related to formal arithmetic. In particular, we show that the free GLP-algebra enjoys (an analogue of) the reduction property and some topological algebras satisfy some generalizations of it. The first of these results answers a question suggested by Joost Joosten (private communication).

\section{Reduction property for GLP-algebras}

Recall that theories extending $T$ are associated with the filters of the Lindenbaum algebra $\cL_T$. A theory $U$ is called $\Pi_n^0$-conservative over a theory $V$ (denoted $V\vdash_{\Pi_n^0} U$) if $U\vdash\pi$ implies $V\vdash\pi$, for all $\Pi_n^0$-sentences $\pi$. This defines a transitive reflexive relation on the set of filters of $\cL_T$. We use the same notation for arbitrary subsets $U,V$ of $\cL_T$ when we mean the same relation for the filters generated by these subsets.

Reduction property states that $\Pi^0_{n+1}$-consequences of an element of $\cM_T$ the form $\la n+1\ra\phi$ (of arithmetical complexity $\Pi^0_{n+2}$) can be axiomatized by a sequence of iterated reflection principles $\{Q_n^k(\phi):k<\gw\}$, where
$$Q^0_n(\phi)= \top, \quad Q_n^k(\phi)=\la n\ra(\phi\land Q_n^k(\phi)).$$
Intuitively, it means that $\la n+1\ra\phi$ is as weak as possible relative to $\la n\ra\phi$ given the constraints of GLP-axioms.

Clearly, for any $\phi$, $n$, and $k$, the elements $Q^k_n(\phi)$ correspond to arithmetical $\Pi_{n+1}^0$-sentences and follow from $\la n+1\ra\phi$. The opposite implication only holds in the sense of  $\Pi^0_{n+1}$-conservativity and under certain assumptions on $T$.

\bt[reduction property, \cite{Bek05}] \label{rp} Suppose $T$ is axiomatized over $\EA$ by an r.e.\ set of $\Pi^0_{n+1}$-sentences. Then
  $$\{Q_n^k(\phi):k<\gw\}\vdash_{\Pi_{n+1}^0} \la n+1\ra\phi.$$
\et

In order to generalize this property to arbitrary GLP-algebras we first introduce some useful notation.
For $A\subseteq \cM$ define $A\vdash x$, if there are $a_1,\dots, a_k\in A$ such that $a_1\land\dots \land a_k\leq x$, in other words, if $x$ belongs to the filter generated by $A$. Define $A\vdash B$ if $A\vdash x$, for all $x\in B$.

For any sets $A,B\subseteq \cM$ define the \emph{$\la n\ra$-conservativity relation} $A\vdash_n B$ by: $$\al{z\in\cM}(B\vdash \la n \ra z \Imp A\vdash \la n\ra z).$$
We put $A\equiv_n B$ iff $A\vdash_n B$ and $B\vdash_n A$.

\bd \label{rpa} A GLP-algebra $\cM$ enjoys the \emph{$\la n \ra$-reduction property} if, for all $x\in \cM$,
$$\{Q_n^k(x):k<\gw\}\vdash_{n} \la n+1\ra x.$$
\ed

Since the elements $\la n\ra z$ of the algebra $\cM_T$, for any $z\in\cM_T$, have complexity $\Pi^0_{n+1}$, Theorem~\ref{rp} shows that $\cM_T$ enjoys the $\la n\ra$-reduction property, for all $n\geq m$, provided $T$ is a $\Pi_{m+2}^0$-axiomatized extension of $\EA$. On the other hand, the following observation shows that the converse also holds.
\bl
If $\cM_T$ satisfies the $\la n\ra$-reduction property then, for any $\phi\in\cM_T$, $$\{Q_n^k(\phi):k<\gw\}\vdash_{\Pi_{n+1}^0} \la n+1\ra\phi.$$
\el

\bp\ Assume $\pi\in\Pi_{n+1}^0$, $\phi\in\cM_T$ and $\la n+1\ra\phi\vdash \pi$. Since $\pi\in\Pi_{n+1}^0$ we have $\la n+1\ra\top\vdash \pi\to \la n\ra\pi$. Hence $\la n+1\ra \phi\vdash \la n\ra \pi$. By the $\la n\ra$-reduction property we infer $\{Q^k_k(\phi):k<\gw\}\vdash \la n\ra\pi\vdash \pi$, as required. \ep

Thus, the abstract reduction property stated in Definition \ref{rpa} is equivalent to the original one for the class of reflection algebras.

The $\la n\ra$-reduction property for a GLP-algebra $\cM$ implies an apparently stronger conservation result.
Let $\Pi_{n+1}(\cM)$ denote the closure under $\lor$, $\land$ of the following subset of $\cM$:
$$\{\top,\bot\}\cup\{\la k \ra z:k\leq n, z\in \cM\}\cup\{[k]z:k<n, z\in \cM\}.$$
Notice that these elements always represent $\Pi_{n+1}^0$-sentences in $\cM_T$. In general, it is not true that $\Pi_{n+1}(\cM_T)$ coincides with the set of all equivalence classes of $\Pi_{n+1}^0$-sentences in $\cM_T$.
By the so-called Friedman--Goldfarb--Harrington (FGH) principle, any $\Pi_{1}^0$-sentence below $\la 0\ra\top$ is equivalent to a sentence of the form $\la 0\ra\phi$, for some $\phi$. It follows that $\Pi_1(\cM_T)$ consists of the equivalence classes of $\Pi_1^0$-sentences implying the consistency assertion for $T$. If $T$ contains the collection schema $B\Sigma_{n}$, a suitable generalization of the FGH-principle holds for $\la n\ra$ in $\cM_T$ (see \cite{Joo15}).

Let $A,B\subseteq \cM$, define $A\vdash_{\Pi_{n+1}(\cM)} B$ if $\al{z\in\Pi_{n+1}(\cM)}(B\vdash z \Imp A\vdash z).$
\bt \label{pi}
Suppose $\cM$ enjoys the $\la n\ra$-reduction property, then
for any $\phi\in\cM$, $$\{Q_n^k(\phi):k<\gw\}\vdash_{\Pi_{n+1}(\cM)} \la n+1\ra\phi.$$
\et

We omit the proof. Theorem \ref{pi} explains why we have not chosen some apparently larger class of modal formulas to represent $\Pi_{n+1}^0$-sentences in Definition~\ref{rpa}. Thus, it appears that $\la n\ra$-reduction property is the right analogue of the reduction property for the reflection algebras in arithmetic, even though the notion of $\Pi_{n+1}^0$-conservativity could be stronger than conservativity for the class of all sentences of the form $\la n\ra \phi$.

\section{Reduction property for free GLP-algebras}

In this section we assume the familiarity with some notions from \cite{Bek10,Bek11} and stick to logical rather than algebraic notation.  We work with Japaridze's logic $\GLP$ and with its fragment $\J$ which has nice Kripke semantics and to which $\GLP$ is reducible by Theorem 4 of \cite{Bek10}.

Our goal is the following theorem.

\bt The free GLP-algebra on any number of generators enjoys the $\la n\ra$-reduction property, for all $n$.
\et

Let $dp(\phi)$ denote the modal depth of a GLP-formula $\phi$. Let $\sim_n$ denote the $n$-bisimilarity equivalence relation on a given model. Recall that $\sim_n$ respects the forcing of formulas of modal depth $\leq n$, and the number of equivalence classes of $\sim_n$ on any model is bounded by a function of $n$ and the number of variables considered.

The following lemma yields a proof of the theorem.
\bl
Suppose $\GLP\vdash \la m+1\ra\psi \to \la m \ra\phi$. Then there is a $k$ such that $\GLP\vdash Q_m^k(\psi)\to \la m\ra\phi$. Moreover, the bound $k$ only depends on
$d=\max(dp(\phi),dp(\psi))$ and the number of variables in $\phi$ and $\psi$.
\el

\bp\ Select a $k$ larger than the number of equivalence classes of $\sim_d$. Assume $\GLP\nvdash Q_m^k(\psi)\to \la m\ra\phi$. Let, as in \cite{Bek10} or \cite{Bek11}, $M(A)$ denote the conjunction of instances of the monotonicity schema $[i]\theta\to [j]\theta$, for all subformulas $[i]\theta$ of $A$ and all $j$ such that $r\geq j>i$, where $r$ is the maximal modality number occurring in $A$. Further, let $M^+(A):=M(A)\land \bigwedge_{i\leq r} [i]M(A)$. Obviously, for any formula $A$, $M^+(A)$ is a theorem of $\GLP$. Hence,
$$\J\nvdash M^+(\la m+1\ra\psi \to \la m \ra\phi)\land Q_m^k(\psi)\to \la m\ra\phi.$$

Let $\cW$ be a rooted J-model of $M^+(\la m+1\ra\psi \to \la m \ra\phi)\land Q_m^k(\psi)\land \neg\la m\ra\phi$. Let $0$ denote its hereditary root. Since $\cW,0\fc Q_m^k(\psi)$, there is a sequence of nodes
$0 R_m a_{k-1} R_m \dots R_m a_0$ such that $\cW,a_i\fc \psi\land Q_m^i(\psi)$, for each $i<k$. By the pigeonhole principle there are $i>j$ such that $a_i\sim_d a_j$.

We denote $a:=a_i$ and $a':=a_j$. Let $\cW_a$ denote the submodel of $\cW$ generated by $a$, and let $\ga$ denote the $(m+1)$-plane generated by $a$. W.l.o.g.\ we may assume that $a$ is the hereditary root of $\ga$, thus $$\ga:=\{x\in \cW: \ex{s>m} a R_{s} x\}\cup\{a\}.$$ Similarly, let $$\ga':=\{x\in \cW: \ex{s>m} a' R_{s} x\}\cup\{a'\}.$$ denote the $(m+1)$-plane generated by $a'$. Further, let $\cW'$ be the minimal J-model obtained from $\cW_a$ by adding a new node $b$ such that $b R_{m+1} x$ for all $x\in \ga$. In particular, $bR_i x$ iff $aR_i x$, for all $i\leq m$, and for no $x$ do we have $b R_j x$ if $j>m+1$. Thus, $b$ is the new hereditary root of $\cW'$.

We claim that $$\cW',b\nfc M^+(\la m+1\ra\psi \to \la m \ra\phi)\land \la m+1\ra\psi\to \la m\ra\phi.$$ This yields $\GLP\nvdash \la m+1\ra\psi\to \la m\ra\phi$, as required.

Clearly, the forcing in $\cW_a$ is the same as in the corresponding part of $\cW$ and of $\cW'$, hence we have
$\cW_a\models M^+(\la m+1\ra\psi \to \la m \ra\phi)$, and $\cW_a, x\models \neg\phi$, for all $x$ such that $bR_m x$ or $b R_{m+1} x$, whence $\cW',b\fc [m]\neg\phi$. Obviously, $\cW',b\fc\la m+1\ra\psi$. This, it remains for us to show that $\cW',b\fc M(\phi)$, which is the content of the following lemma.

\bl \label{emm} For any subformula $[i]\theta$ of the formula $\la m+1\ra\psi\to \la m\ra\phi$ and any $j>i$ there holds $\cW',b\fc [i]\theta\to [j]\theta$.
\el

\bp\ Recall that $\la i\ra\eta$ abbreviates $\neg[i]\neg\eta$, for any formula $\eta$. We consider the following cases.

If $[i]\theta= [m]\neg\phi$ then the only interesting case is $j=m+1$, but we already know that $\cW',b\fc [m+1]\neg\phi$.

If $[i]\theta= [m+1]\neg\psi$ then we know that $\cW', b\nfc [m+1]\neg\psi$.

Suppose $[i]\theta$ is a subformula of $\phi$ or $\psi$. We consider the following subcases.

\ben
\item If $j>m+1$ then trivially $\cW', b\fc [j]\theta$.

\item If $j\leq m$ we obtain: $\cW',b\fc [i]\theta$ implies $\cW',a\fc [i]\theta$, hence $\cW,a\fc [i]\theta$ and $\cW,a\fc [j]\theta$, because
    $\cW,a\fc [i]\theta\to [j]\theta$. Since $i<j\leq m$ this yields $\cW',b\fc[j]\theta$.

\item If $j=m+1$ and $i<m$ then $\cW',b\fc [i]\theta$ implies $\cW,0\fc [i]\theta$, hence $\cW,0\fc [m]\theta$. It follows that $\cW,x\fc \theta$, for each $x\in\ga$, and therefore $\cW',b\fc [m+1]\theta$.

\item If $j=m+1$ and $i=m$ we first notice that $dp(\theta)\leq d-1$ and, since $a\sim_d a'$, there holds
$$\al{x\in\ga}\ex{y\in\ga'} x\sim_{d-1} y.$$
Therefore, $\cW',b\fc [m]\theta$  implies $\cW,y\fc \theta$, for each $y\in\ga'$. Hence, $\cW,x\fc \theta$, for each $x\in\ga$. Thus, $\cW',b\fc [m+1]\theta$.

\een
This proves Lemma \ref{emm}. \ep

Lemma \ref{emm} together with the previous shows that
$$\J\nvdash M^+(\la m+1\ra\psi \to \la m \ra\phi)\land \la m+1\ra\psi\to \la m\ra\phi.$$ Hence, by Theorem 4 of \cite{Bek10}, $\GLP\nvdash \la m+1\ra\psi\to \la m\ra\phi.$
\ep

\section{Reduction property for GLP-spaces}
The following definition comes from \cite{BBI09}.
\bd A \emph{GLP-space} is a nonempty set equipped with a sequence of topologies  $(X,\{\tau_n:n<\gw\})$ such that $(\cP(X),\{d_n:n<\gw\})$ is a GLP-algebra.
\ed

Here, $d_n$ denotes the topological derivative operator w.r.t.\ topology $\tau_n$, that is, $d_n(A)$ is the set of all limit points of a subset $A$ in $X$. Similarly, $c_n$ will denote the closure operator w.r.t.\ $\tau_n$.

It is well-known that in a GLP-space, for all $n<\gw$,
\ben \item $\tau_n$ is scattered, that is, every non-empty subspace has a isolated point;
\item $\tau_n\subseteq\tau_{n+1}$;
\item For each $A\subseteq X$, $d_n(A)$ is $\tau_{n+1}$-open.
\een

The main example of a GLP-space is the \emph{ordinal GLP-space}, that is, the space $(\Omega,\{\tau_n:n<\gw\})$ where $\Omega$ is an ordinal, $\tau_0$ is the left topology and $\tau_{n+1}$ is generated by $\tau_n\cup \{d_n(A):A\subseteq \Omega\}$. We notice that, for each $n<\gw$, $\tau_n$ is zero-dimensional, as the sets $d_n(A)$ are clopen in the next topology $\tau_{n+1}$. Also, each $\tau_n$ for $n\geq 1$ is $T_3$.

We mention without proof the following characterization.
\bpr
Let $(X,\{\tau_n:n<\gw\})$ be a $T_3$ GLP-space. For any $A,B\subseteq X$, $A\vdash_n B$ iff ($c_n(A)\subseteq c_n(B)$ or $B\nsubseteq d_n(X)$).
\epr

\renewcommand{\DD}{d}
Topological analogs of the terms $Q_n^k(\phi)$ are defined as follows, where we generalize to transfinite iterations. Let $A\subseteq X$, $\ga$ an ordinal, $\gl$ a limit ordinal; define:
$$d_n^0[A]=X;\quad d_n^{\ga+1}[A]=d_n(d_n^\ga[A]\cap A);\quad
d_n^\gl[A]=\bigcap_{\ga<\gl} d_n^\ga[A].$$ We also note that $d_n^\ga[X]$ is the familiar Cantor--Bendixson sequence for $(X,\tau_n)$. If $(X,\tau_n)$ is scattered the sequence $d_n^\ga[A]$ is a strictly decreasing sequence of $\tau_n$-closed sets, hence we have $d_n^\ga[A]=\emptyset$, for some $\ga$.

Unwinding the definitions we see that a GLP-algebra $(\cP(X),\{d_k:k<\gw\})$ satisfies the $\la n\ra$-reduction property iff
$$\{d_n^k[A]:k<\gw\}\vdash_n d_{n+1}(A).$$ This is, in general, stronger than saying $d_n^\gw[A]\vdash_n d_{n+1}(A).$
(Consider, for example, the left topology and the interval topology on $\gw$.)
However, if $(X,\tau_{n+1})$ is compact then, for all $A,B\subseteq X$ and limit ordinals $\gl$,
 $$d_n^\gl[A]\subseteq d_n(B) \iff \ex{\ga<\gl} d_n^{\ga}[A]\subseteq d_n(B).$$ Indeed, $X\setminus d_n(B)$ is $\tau_{n+1}$-closed, hence compact, and the left hand side means that $X\setminus d_n^\gl[A]=\bigcup_{\ga<\gl} (X\setminus d_n^\ga[A])$ is its open cover.

 The following more general definition seems to be working well also in the non-compact case.

\bd A GLP-space $X$ satisfies (weak)
\emph{$\ga$-reduction property for $d_n$} if, for each subset $A\subseteq X$,
$$ d_n^\ga[A]\vdash_n d_{n+1}(A).$$
\ed
If $(X,\tau_n)$ is $T_3$, the weak $\ga$-reduction property for $d_n$
is equivalent to the identity $c_n(d_{n+1}(A))= d_n^\ga[A],$
for any $A\subseteq X$.

\bt The ordinal GLP-space $(\Omega,\{\tau_k:k<\gw\})$ (for $\Omega$ sufficiently large) satisfies
\benr
\item weak $\gw$-reduction property for $d_0$;
\item weak $\gw_1$-reduction property for $d_1$.
\eenr
\et

We conjecture that more generally weak $\kappa$-reduction property for $d_n$ holds, where $\kappa$ is the first limit point of $(\Omega,\tau_{n+1})$ (if such a point exists). The existence of limit points for $\tau_n$ for $n>2$ is a large cardinal hypothesis independent of ZFC.

\bibliographystyle{plain}

\input{reduction-en.bb}
\end{document}

%% file: macros.tex
\newtheorem{theorem}{Theorem}
\newtheorem{lemma}{Lemma}

\newtheorem{corollary}[lemma]{Corollary}
\newtheorem{proposition}[lemma]{Proposition}

\newtheorem{definition}{Definition}
\newtheorem{example}[lemma]{Example}
\newtheorem{remark}[lemma]{Remark}

\newcommand{\bl}{\begin{lemma}}
\newcommand{\el}{\end{lemma}}
\newcommand{\bt}{\begin{theorem}}
\newcommand{\et}{\end{theorem}}
\newcommand{\bcor}{\begin{corollary}}
\newcommand{\ecor}{\end{corollary}}
\newcommand{\bp}{\proof{.}}
\newcommand{\ep}{\eop}
\newcommand{\bpr}{\begin{proposition}}
\newcommand{\epr}{\end{proposition}}
\newcommand{\brem}{\begin{remark} \em}
\newcommand{\erem}{\end{remark}}
\newcommand{\bd}{\begin{definition} \em}
\newcommand{\ed}{\end{definition}}
\newcommand{\bex}{\begin{example} \em
}
\newcommand{\eex}{\end{example}}
\newcommand{\beq}{\begin{equation} }
\newcommand{\eeq}{\end{equation}}

\newcommand{\bi}{\begin{itemize}
  }
\newcommand{\ei}{\end{itemize}}
\newcommand{\ben}{\begin{enumerate} }
\newcommand{\een}{\end{enumerate} }

\newcommand{\refeq}[1]{(\ref{#1})}
\newenvironment{enumr}{

\begin{enumerate}     }{\end{enumerate}

}

\newcommand{\benr}{\begin{enumr}
  }
\newcommand{\eenr}{
\end{enumr}}

\newcommand{\al}[1]{\forall #1\:}
\newcommand{\ex}[1]{\exists #1\:}

\newlength{\hilflh}

\renewcommand{\emptyset}{\varnothing}

\newcommand{\cL}{{\mathcal L}}

\newcommand{\cP}{{\mathcal P}}

\newcommand{\cM}{{\mathcal M}}

\newcommand{\cW}{{\mathcal W}}

\newcommand{\ga}{\alpha}

\renewcommand{\ge}{\varepsilon}
\newcommand{\gl}{\lambda}

\newcommand{\gw}{\omega}

\renewcommand{\phi}{\varphi}

\newcommand{\GLP}{\mathbf{GLP}}

\newcommand{\DD}{\mbox{\bf D}}

\newcommand{\PA}{\mathsf{PA}}
\newcommand{\EA}{\mathrm{EA}}

\newcommand{\la}{\langle}
\newcommand{\ra}{\rangle}

\newcommand{\fc}{\Vdash}      
\renewcommand{\models}{\vDash}      
\newcommand{\nfc}{\nVdash}

\newcommand{\Imp}{\Rightarrow}

\newcommand{\J}{\mathbf{J}}

\renewcommand{\leq}{\leqslant}
\renewcommand{\geq}{\geqslant}

\newcommand{\eop}{$\boxtimes$ \protect\par \addvspace{\topsep}}
\newcommand{\proof}[1]{\protect\par\addvspace{\topsep}\noindent {\bf Proof#1}}

